\documentclass[12pt,a4paper]{article} 
\usepackage[T1]{fontenc}
\usepackage{mathptmx}
\usepackage[dvips]{graphicx}
\usepackage{psfrag}
\usepackage{amsmath,amssymb}        % moduli per le estensioni ed i fonts AMS
\usepackage{latexsym}
\usepackage{url}
\usepackage{booktabs}

\def\RR{\hbox{I\kern-.2em\hbox{R}}}

\newcommand{\eqnsection}{
   \renewcommand{\theequation}{{\thesection.\arabic{equation}}}
   \makeatletter
   \csname @addtoreset\endcsname{equation}{section}
   \makeatother}
\eqnsection
\title{The Iterative Transformation Method for \\
the Sakiadis Problem}
\author{Riccardo Fazio \\
Department of Mathematics and Computer Science \\
University of Messina \\
Viale F. Stagno D'Alcontres 31, 98166 Messina, Italy\\
{\small e-mail: rfazio@unime.it \ \ \ 
home-page: http://mat521.unime.it/fazio}}
\date{October 8, 2013 and in revised form \today}
\begin{document}               
\maketitle

\begin{abstract}
In a transformation method the numerical solution of a given boundary value problem is obtained by solving one or more related initial value problems.
This paper is concerned with the application of the iterative transformation method to the Sakiadis problem. 
This method is an extension of the T{\"o}pfer's non-iterative algorithm developed as a simple way to solve the celebrated Blasius problem.
As shown by this author [Appl. Anal., {\bf 66} (1997) pp. 89-100] the method provides a simple numerical test for the existence and uniqueness of solutions.
Here we show how the method can be applied to problems with a homogeneous boundary conditions at infinity and in particular we solve the Sakiadis problem of boundary layer theory. 
Moreover, we show how to couple our method with Newton's root-finder.
The obtained numerical results compare well with those available in literature.
The main aim here is that any method developed for the Blasius, or the Sakiadis, problem might be extended to more challenging or interesting problems.
In this context, the iterative transformation method has been recently applied to compute the normal and reverse flow solutions of Stewartson for the Falkner-Skan model [Comput. \& Fluids, {\bf 73} (2013) pp. 202-209].
\end{abstract}
% \centerline{Abbreviated title (for running head): Numerical solution
% of moving boundary parabolic problems.}



\section{Introduction}\label{Introduction}
In a transformation method the numerical solution of a given boundary value problem is obtained by solving one or more related initial value problems (IVPs).
In this context the classical example is the Blasius problem of boundary layer theory.  
In the Blasius problem the governing differential equation and the two initial conditions are invariant under the scaling group of transformations
\begin{equation}\label{eq:scalinga}
f^* = \lambda^{-\alpha} f \ , \qquad \eta^* = \lambda^{\alpha} \eta \ ,   
\end{equation}
where $\lambda $ is the group parameter and $\alpha \ne 0$.
Moreover, the non-homogeneous asymptotic boundary condition is not invariant with respect to (\ref{eq:scalinga}). 
This kind of invariance was used by T\"opfer \cite{Topfer:1912:BAB} to define a non-iterative transformation method (ITM) for the Blasius problem by transforming the boundary conditions to initial conditions and rescaling the obtained numerical solution. 

This paper is concerned with the application an ITM to the Sakiadis problem.
The main aim here is that any method developed for the Blasius, or the Sakiadis, problem might be extended to more challenging or interesting problems.
In this context, the iterative transformation method has been recently applied to compute the normal and reverse flow solutions of Stewartson \cite{Stewartson:1954:FSF,Stewartson:1964:TLB} for the Falkner-Skan model \cite{Fazio:2013:BPF}.

The Sakiadis problem is a variant of Blasius problem that cannot be solved by a non-ITM.
In fact, one of the initial conditions is not invariant and the asymptotic boundary condition, being homogeneous, is invariant with respect to the scaling transformations (\ref{eq:scalinga}).
Therefore, as noted by Na \cite[pp. 160-164]{Na:1979:CME}, it is not possible to rescale an initial value solution to the given asymptotic boundary condition. 
Moreover, the non-ITM cannot be applied when the governing differential equation is not invariant with respect to a scaling group of point transformations.
To overcome this drawback the ITM was defined in \cite{Fazio:1994:FSE,Fazio:1996:NAN} for the numerical solution of the Falkner-Skan model and of other problems in boundary layer theory.

Here we show how the ITM can be applied to problems with a homogeneous boundary conditions at infinity.
Moreover, we show how to couple our method with the Newton's root-finder.
This ITM has been applied to several problems of interest: 
free boundary problems \cite{Fazio:1991:ITM,Fazio:1997:NTE},
a moving boundary hyperbolic problem \cite{Fazio:1992:MBH},
the Falkner-Skan equation in \cite{Fazio:1994:FSE,Fazio:1996:NAN,Fazio:2013:BPF}, one-dimensional parabolic moving boundary problems \cite{Fazio:2001:ITM,Fazio:2013:SII}, two variants of the Blasius problem \cite{Fazio:2009:NTM}, namely: a boundary layer problem over moving plates, studied first by Klemp and Acrivos \cite{Klemp:1972:MBL}, and a boundary layer problem with slip boundary condition, that has found application to the study of gas and liquid flows at the micro-scale regime \cite{Gad-el-Hak:1999:FMM,Martin:2001:BBL}, a parabolic problem on unbounded domain \cite{Fazio:2010:MBF}.
Furthermore, as shown in \cite{Fazio:1997:NTE}, the ITM provides a simple numerical test for the existence and uniqueness of solutions.

\section{Blasius and Sakiadis problems}\label{Sakiadis}
Within boundary-layer theory, the model describing the steady plane flow of a fluid past a thin plate, is given by
\begin{align}\label{PDE-model}
& {\displaystyle \frac{\partial u}{\partial x}} +
{\displaystyle \frac{\partial v}{\partial y}} = 0 \ ,  \nonumber \\[-1.2ex]
& \\[-1.2ex]
& u {\displaystyle \frac{\partial u}{\partial x}} +
v {\displaystyle \frac{\partial u}{\partial y}} = \nu
{\displaystyle \frac{\partial^2 u}{\partial y^2}} \ ,  \nonumber 
\end{align}
where the governing differential equations, namely conservation of mass and momentum, are the 
steady-state 2D Navier-Stokes equations under the boundary layer approximations: $ u \gg v $ and the flow has a very thin layer attached to the plate,
$ u $ and $ v $ are the velocity components of the fluid in the
$ x $ and $ y $ direction, and $ \nu $ is the viscosity of the fluid.
The boundary conditions for the velocity field are
\begin{align}\label{PDE_BCs-Blasius}
& u(x, 0) = v(x, 0) = 0 \ , \qquad u(0, y) = U_\infty  \ , \nonumber \\[-1.2ex]
&\\[-1.2ex]
& u(x, y) \rightarrow U_{\infty} \quad \mbox{as}
\quad y \rightarrow \infty \ , \nonumber 
\end{align}
for the Blasius flat plate flow problem \cite{Blasius:1908:GFK}, where $ U_{\infty} $ is the main-stream velocity, and 
\begin{align}\label{PDE_BCs-Sakiadis}
& u(x, 0) = U_p \ , \qquad  v(x, 0) = 0 \ , \nonumber \\[-1.2ex]
&\\[-1.2ex]
& u(x, y) \rightarrow 0 \quad \mbox{as}
\quad y \rightarrow \infty \ , \nonumber 
\end{align}
for the classical Sakiadis flat plate flow problem \cite{Sakiadis:1961:BLBa,Sakiadis:1961:BLBb}, where $ U_p $ is the plate velocity, respectively.
The boundary conditions at $ y = 0 $ are based on the assumption that neither slip nor mass
transfer are permitted at the plate whereas the remaining boundary condition means that the velocity $ v $ tends to the main-stream velocity $ U_{\infty} $ asymptotically or gives the prescribed velocity of the plate $U_p$.

Introducing a similarity variable $\eta$ and a dimensionless stream function $f(\eta)$ as
\begin{equation}\label{eq:simvar}
\eta = y \sqrt{\frac{U}{\nu x}} \ , \qquad u = U \frac{df}{d\eta} \ , \qquad v = \frac{1}{2}\sqrt{\frac{U \nu}{x}}\left(\eta \frac{df}{d\eta}-f\right) \ ,
\end{equation}
we have
\begin{equation}\label{eq:simvar1}
\frac{\partial u}{\partial x} = -\frac{U}{2}\frac{\eta}{x}\frac{d^2f}{d\eta^2} \ , \qquad \frac{\partial v}{\partial y} = \frac{U}{2}\frac{\eta}{x}\frac{d^2f}{d\eta^2}
\end{equation}
and the equation of continuity, the first equation in (\ref{PDE-model}), is satisfied identically.

On the other hand, we get
\begin{equation}\label{eq:simvar2}
\frac{\partial u}{\partial y} = U\frac{d^2f}{d\eta^2} \sqrt{\frac{U}{\nu x}} \ , \qquad \frac{\partial^2 u}{\partial y^2} = \frac{U^2}{\nu x}\frac{d^3f}{d\eta^3} \ .
\end{equation}
Let us notice that, in the above equations $U = U_\infty$ represents Blasius flow, whereas $U = U_p$ indicates Sakiadis flow, respectively.

By inserting these expressions into the momentum equation, the second equation in (\ref{PDE-model}), we get 
\begin{equation}\label{eq:Blasius}
{\displaystyle \frac{d^3 f}{d \eta^3}} + \frac{1}{2} f
{\displaystyle \frac{d^{2}f}{d\eta^2}} = 0 \ ,
\end{equation}
to be considered along with the transformed boundary conditions
\begin{equation}\label{eq:Blasius:BCs}
f(0) = {\displaystyle \frac{df}{d\eta}}(0) = 0 \ , \qquad
{\displaystyle \frac{df}{d\eta}}(\eta) \rightarrow 1 \quad \mbox{as}
\quad \eta \rightarrow \infty \ , \nonumber
\end{equation}
for the Blasius flow, and
\begin{equation}\label{eq:Sakiadis:BCs}
f(0) = 0 \ , \qquad {\displaystyle \frac{df}{d\eta}}(0) = 1 \ , \qquad
{\displaystyle \frac{df}{d\eta}}(\eta) \rightarrow 0 \quad \mbox{as}
\quad \eta \rightarrow \infty \ , \nonumber 
\end{equation}
for the Sakiadis flow, respectively.

Blasius main interest was to compute the value of 
the velocity gradient at the plate (the wall shear or skin friction coefficient):
\[
 \lambda = \frac{d^2f}{d\eta^2}(0) \ .
\]
To compute this value, Blasius used a formal series solution around $\eta=0$ and an asymptotic expansions for large values of $\eta$, adjusting the constant $\lambda$ so as to connect both expansions in a middle region.
In this way, Blasius obtained the (erroneous) bounds $ 0.3315 < \lambda < 0.33175$.

A few years later, T{\"o}pfer \cite{Topfer:1912:BAB} revised the work by Blasius and solved numerically the Blasius problem, using a non-ITM.
He then arrived, without detailing his computations, at the value $\lambda \approx 0.33206$, contradicting the bounds reported by Blasius.  

Indeed, T{\"o}pfer solved the IVP for the Blasius equation once.
At large but finite $ \eta_j^* $, ordered so that $ \eta_j^* < \eta_{j+1}^* $, he computed the corresponding scaling parameter $ \lambda_j $.  
If two subsequent values of $ \lambda_j $ agree within a specified accuracy, then $ \lambda $ is approximately equal to the common value of the $ \lambda_j $, otherwise, he marched to a larger value of $ \eta^* $ and tried again.
Using the classical fourth order Runge-Kutta method, as given by Butcher \cite[p. 166]{Butcher}, and a grid step $ \Delta \eta^* = 0.1$ T{\"o}pfer was able, only by hand computations, to determine $ \lambda $ with an error less than $ 10^{-5} $. 
To this end he used the two {\it truncated boundaries} $\eta_1^* = 4$ and $\eta_2^* = 6$.
For the sake of simplicity we follow T\"opfer and apply some preliminary computational tests to find a suitable value for the truncated boundary.

Sakiadis studied the behaviour of boundary layer flow, due to a moving flat plate immersed in an otherwise quiescent fluid, \cite{Sakiadis:1961:BLBa,Sakiadis:1961:BLBb}.
He found that the wall shear is about 34\% higher for the Sakiadis flow compared to the Blasius case. 
Later, Tsou and Goldstein \cite{Tsou:1967:FHT} made an experimental and theoretical treatment of Sakiadis problem to prove that such a flow is physically realizable.

\section{Extension of T{\"o}pfer algorithm: the ITM}
Within this section we explain how it is possible to extend T{\"o}pfer algorithm to the Sakiadis problem, that we rewrite here for the reader convenience
\begin{eqnarray}\label{eq:Sakiadis}
&{\displaystyle \frac{d^3 f}{d \eta^3}} + \frac{1}{2} f
{\displaystyle \frac{d^{2}f}{d\eta^2}} = 0 \nonumber \\[-1.5ex]
&\\[-1.5ex]
&f(0) = 0 \ , \qquad {\displaystyle \frac{df}{d\eta}}(0) = 1 \ , \qquad
{\displaystyle \frac{df}{d\eta}}(\eta) \rightarrow 0 \quad \mbox{as}
\quad \eta \rightarrow \infty \ . \nonumber 
\end{eqnarray}

In order to define the ITM we introduce the extended problem
\begin{eqnarray}\label{eq:Sakiadis:mod}
&{\displaystyle \frac{d^3 f}{d \eta^3}} + \frac{1}{2} f
{\displaystyle \frac{d^{2}f}{d\eta^2}} = 0 \nonumber \\[-1.5ex]
& \\[-1.5ex]
&f(0) = 0 \ , \qquad {\displaystyle \frac{df}{d\eta}}(0) = h^{1/2} \ , \qquad
{\displaystyle \frac{df}{d\eta}}(\eta) \rightarrow 1 - h^{1/2} \quad \mbox{as}
\quad \eta \rightarrow \infty \ . \nonumber
\end{eqnarray}
In (\ref{eq:Sakiadis:mod}), the governing differential equation and the two initial conditions are invariant, the asymptotic boundary condition is not invariant, with respect to the extended scaling group
\begin{equation}\label{eq:scaling}
f^* = \lambda f \ , \qquad \eta^* = \lambda^{-1} \eta \ , \qquad 
h^* = \lambda^{4} h \ .   
\end{equation}
Moreover, it is worth noticing that the extended problem (\ref{eq:Sakiadis:mod}) reduces to the Sakiadis problem  (\ref{eq:Sakiadis}) for $h=1$.
This means that in order to find a solution of the Sakiadis problem we have to find a zero of the so-called {\it transformation function} 
\begin{equation}\label{eq:AA2.9}
\Gamma (h^{*}) = \lambda^{-4} h^* - 1 \ , 
\end{equation}  
where the group parameter $ \lambda $ is defined with the formula
\begin{equation}\label{eq:lambda}
\lambda = \left[\displaystyle \frac{df^*}{d\eta^*}(\eta_\infty^*)+{h^*}^{1/2}\right]^{1/2} \ ,
\end{equation}
and to this end we can use a root-finder method.

Let us notice that $\lambda$ and the transformation function are defined implicitly by the solution of the IVP
\begin{eqnarray}\label{eq:Sakiadis:IVP}
&{\displaystyle \frac{d^3 f^*}{d \eta^{*3}}} + \frac{1}{2} f^*
{\displaystyle \frac{d^{2}f^*}{d\eta^{*2}}} = 0 \nonumber \\[-1.5ex]
& \\[-1.5ex]
&f^*(0) = 0 \ , \quad {\displaystyle \frac{df^*}{d\eta^*}}(0) = {h^*}^{1/2}, \quad
{\displaystyle \frac{d^2f^*}{d\eta^{*2}}}(0) = \pm 1  \nonumber \ .
\end{eqnarray}
In particular, we are interested to compute $\frac{df^*}{d\eta^*}(\eta_\infty^*)$, an approximation of the asymptotic value $\frac{df^*}{d\eta^*}(\infty)$, which is used in the definition of $\lambda$ (\ref{eq:lambda}). 

For the ITM we have to follow the steps:
\begin{enumerate}
\item{} we apply a root-finder method to define a sequence $ h^{*}_{j}, $ for $ j = 0, 1, 2, \dots  $ 
Two sequences $\lambda_j$ and $ \Gamma (h^{*}_{j}) $ for 
$ j = 0, 1, 2, \dots, \ $ are defined by equation (\ref{eq:lambda}) and (\ref{eq:AA2.9}), respectively.
\item{} a suitable convergence criterion should be used to verify whether $ \Gamma (h^{*}_{j}) \rightarrow  0 $ as $ j \rightarrow \infty $.
If this is the case, then $\lambda_j$ converges to the correct value of $\lambda$ in the same limit.
\item{} a solution of the original problem can be obtained by rescaling to $ h = 1 $.
In particular, we have that
\begin{equation}\label{eq:MIC}
{\displaystyle \frac{d^2f}{{d\eta}^2}}(0) = \lambda^{- 3}
{\displaystyle \frac{d^2f^*}{{d\eta^*}^2}}(0) \ ,
\end{equation}
where this $\lambda$ is the limit value mentioned in the previous step.
\end{enumerate}

Several questions are of interest.
As far as the missing initial condition is concerned, are we allowed to use the value
\[
\frac{d^2f^*}{{d\eta^*}^2} (0) = 1 \ ,
\]
suggested to T{\"o}pfer, see for instance Na \cite[p. 141]{Na:1979:CME}, by a formal series solution of the Blasius problem?
Indeed, if the first derivative of $f$ is a monotone decreasing function, then the given boundary conditions in (\ref{eq:Sakiadis}) indicate that the second derivative of $f(\eta)$ has to be negative and should go to zero as $\eta$ goes to infinity and this calls for a negative value of the missing initial condition.
Is the solution of the Sakiadis problem (\ref{eq:Sakiadis}) unique?
By studying the behaviour of the transformation function $\Gamma$ we can answer both questions, and this is done in the next section.

\section{Numerical Results}
It is evident that our numerical method is based on the behaviour of the transformation function.
Our interest is to study the behaviour of this function with respect to its independent variable as well as the involved parameters.
We notice that, because of the two terms $h^{1/2}$, which have been introduced in the modified boundary conditions in (\ref{eq:Sakiadis:mod}), we are allowed to consider only positive values of $h^*$. 

\begin{figure}[!hbt]
	\centering
\psfrag{h*}[][]{$ h^* $} 
\psfrag{G}[][]{$ \Gamma(h^*) $} 
\psfrag{df}[][]{${\displaystyle \frac{d^2f^*}{{d\eta^*}^2}} (0) = 1$} 
\includegraphics[width=.8\textwidth]{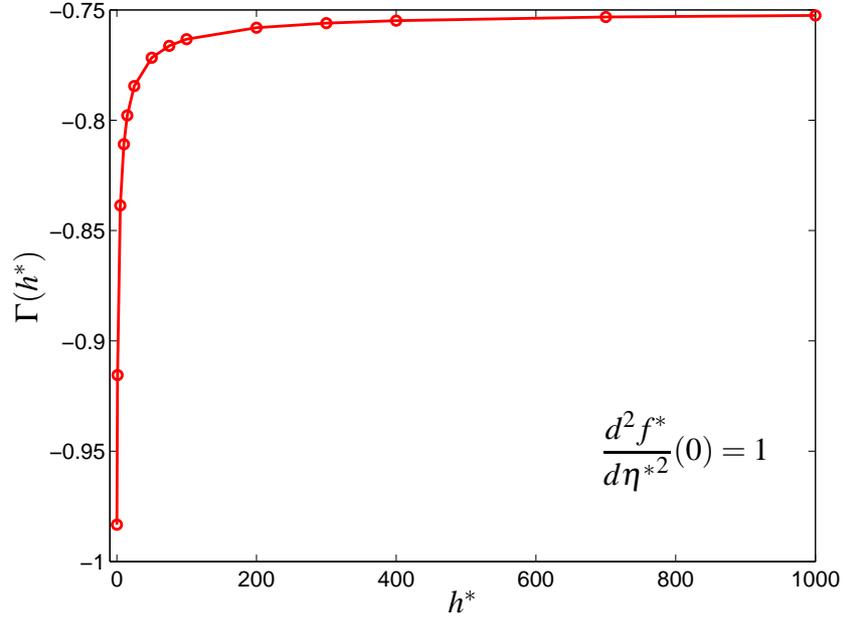} 
\caption{Plot of $\Gamma(h^*)$ for ${\displaystyle \frac{d^2f^*}{{d\eta^*}^2}}(0) = 1$.} 
	\label{fig:Gamma1}
\end{figure}

\begin{figure}[!hbt]
	\centering
\psfrag{h*}[][]{$ h^* $} 
\psfrag{G}[][]{$ \Gamma(h^*) $} 
\psfrag{df}[r][]{${\displaystyle \frac{d^2f^*}{{d\eta^*}^2}} (0) = -1$} 
\includegraphics[width=.9\textwidth]{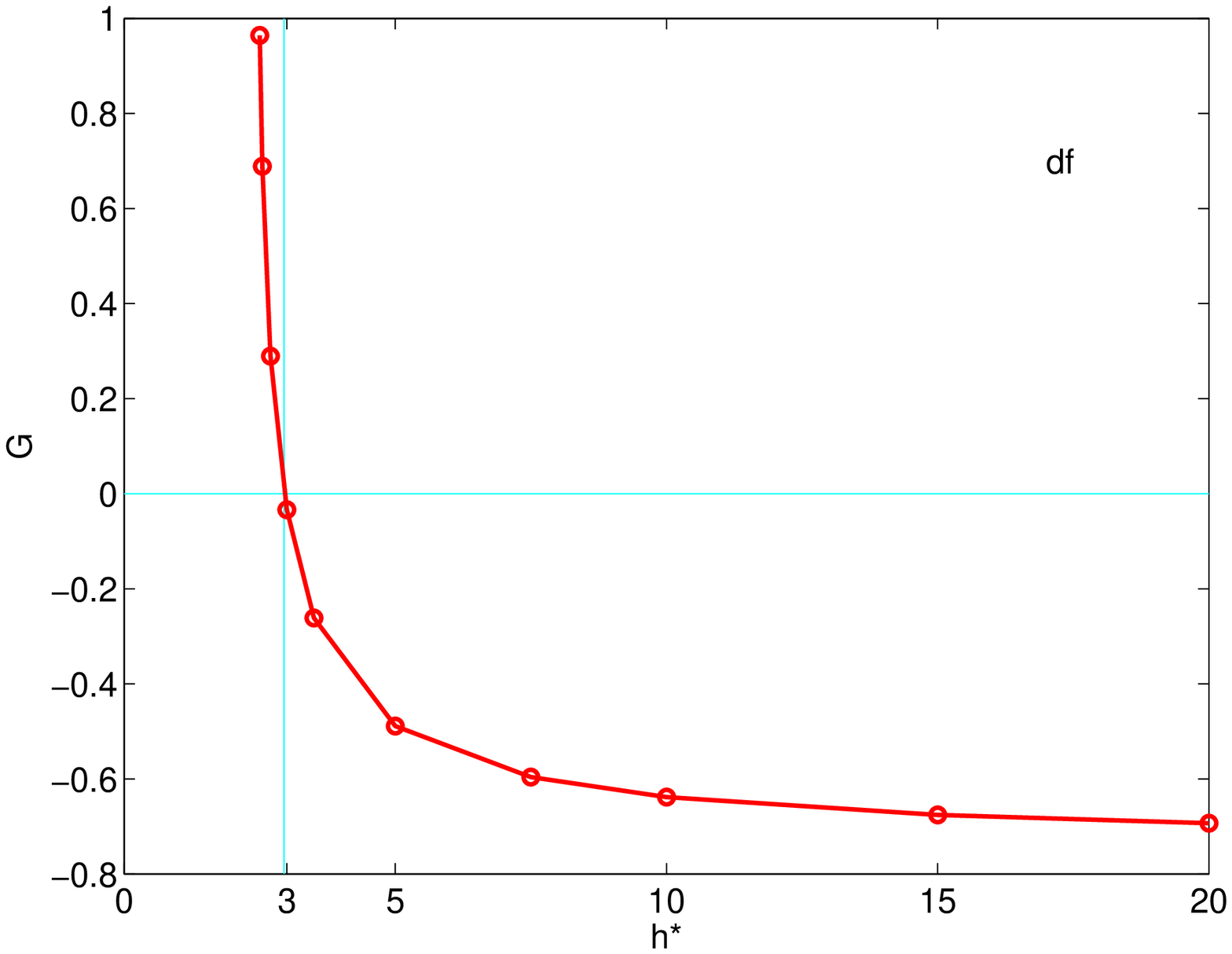} 
\caption{Plot of $\Gamma(h^*)$ for ${\displaystyle \frac{d^2f^*}{{d\eta^*}^2}} (0) = -1$.} 
	\label{fig:Gamma2}
\end{figure}

From our numerical study concerning the dependence of $\Gamma$ with respect to the missing initial condition $\frac{d^2f^*}{{d\eta^*}^2} (0)$ we have used the results plotted on figures \ref{fig:Gamma1}-\ref{fig:Gamma2}.
Each o-symbol represent a numerical solution of the IVP (\ref{eq:Sakiadis:IVP}) with the corresponding value of $h^*$. 
The solid line joining these symbols is used to show the behaviour of the transformation function.  
By considering the results reported on figure \ref{fig:Gamma1} we realize that the missing initial condition cannot be positive.

For a negative missing initial condition the numerical results are shown on figure \ref{fig:Gamma2}.
It is evident from figure \ref{fig:Gamma2} that the transformation function has only one zero and, by a theorem proved in \cite{Fazio:1997:NTE}, this means that the considered problem has one and only one solution.
Moreover, we remark that the tangent to the $\Gamma$ function at its unique zero and the $h^*$ axis define a large angle.
From a numerical viewpoint, this means that the quest for the $h^*$ corresponding to $h=1$ is a well conditioned problem.
 
For a problem, in boundary layer theory, admitting more that one solutions or none, depending on the value of a parameter involved, see \cite{Fazio:2009:NTM} or \cite{Fazio:2013:BPF}.

As far as the numerical results reported in this section are concerned, the ITM was applied by setting the truncated boundary $ \eta_\infty^* = 10$.
Moreover, these results were obtained by an adaptive fourth-order Runge-Kutta IVP solver. 
The adaptive solver uses a relative and an absolute error tolerance, for each component of the numerical solution, both equal to $1 \mbox{D}{-06}$. 
Here and in the following the notation $\mbox{D}-k = 10^{-k}$ means a double precision arithmetic.

\subsection{Secant root-finder}
As a first case the initial value solver was coupled with the simple secant root-finder with a convergence criterion given by
\begin{equation}\label{eq:conv}
\left|\Gamma(h^*)\right| \le 1 \mbox{D}{-09} \ .
\end{equation}

\begin{figure}[!hbt]
	\centering
\psfrag{e}[][]{$ \eta $} 
%\psfrag{c}[][]{$ \frac{df^*}{d\eta^*}(\eta^*), \frac{d^2f^*}{d\eta^{*2}}(\eta^*), \frac{df}{d\eta}(\eta), \frac{d^2f}{d\eta^{2}}(\eta) $} 
\psfrag{f}[][]{$f$} 
\psfrag{df}[][]{$ {\displaystyle \frac{df}{d\eta}} $} 
\psfrag{ddf}[][]{$ {\displaystyle \frac{d^2f}{d\eta^2}} $} 
\includegraphics[width=.9\textwidth]{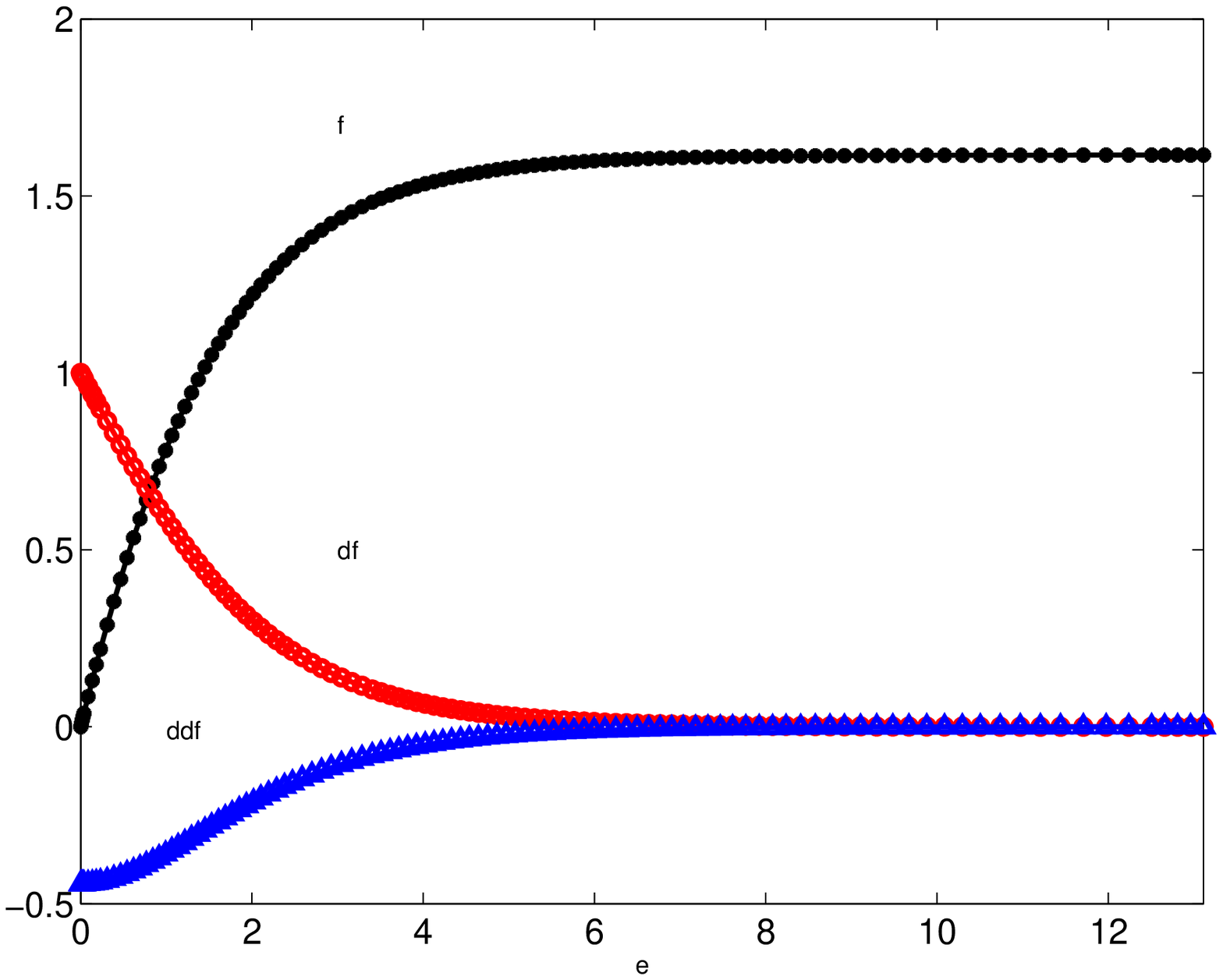} 
\caption{Sakiadis solution via the ITM.} 
	\label{fig:Sakiadis}
\end{figure}

The implementation of the secant method is straightforward.
The only difficulty we have to face is related to the choice of the initial iterates.
In this context the study of the transformation function of figure \ref{fig:Gamma2} can be helpful.
In table \ref{tab:itera} we list the iterations of our ITM.
\begin{table}[!htb]
\caption{Iterations and numerical results: secant root-finder.}
\begin{center}{ \renewcommand\arraystretch{1.3}
%\begin{tabular}{ccccr@{.}lr @{.} lcc} 
\begin{tabular}{rr@{.}lcr@{.}lr@{.}l} 
\toprule%
{$j$}
&
\multicolumn{2}{l}%
{$h_j^*$} 
& {$\lambda_j$}
& \multicolumn{2}{c}%
{$\Gamma(h_j^*)$} 
& \multicolumn{2}{c}%
{${\displaystyle \frac{d^2f}{{d\eta}^2}} (0)$} \\[1.5ex]
\toprule%
% & & & & & \\
$0$  & $2$ & $5$       & $1.061732$ & $0$  & $967343$ & $-0$ & $835517$ \\
$1$  & $3$ & $5$       & $1.475487$ & $-0$ & $261541$ & $-0$ & $311310$ \\
$2$ & $3$ & $287172$ & $1.417981$ & $-0$ & $186906$ & $-0$& $350743$\\
$3$ & $2$ & $754191$ & $1.229206$ & $0$ & $206411$ & $-0$ & $538426$ \\
$4$ & $3$ & $033897$ & $1.339089$ & $-0$ & $056455$ & $-0$ & $416458$ \\
$5$ & $2$ & $973826$ & $1.318081$ & $-0$ & $014749$ & $-0$ & $436690$ \\
$6$ & $2$ & $952581$ & $1.310382$ & $0$ & $001407$ & $-0$ & $444433$ \\
$7$ & $2$ & $954432$ & $1.311058$ & $-3$ & $23\mbox{D}{-05}$ & $-0$ & $443745$ \\
$8$ & $2$ & $954391$ & $1.311043$ & $-6$ & $93\mbox{D}{-08}$ & $-0$ & $443761$ \\
$9$ & $2$ & $954391$ & $1.311043$ & $3$ & $42\mbox{D}{-12}$ & $-0$ & $443761$ \\
\bottomrule
\end{tabular}}
\label{tab:itera}
\end{center}
\end{table}
The last iteration of table \ref{tab:itera} defines our numerical approximation that is shown, for the reader convenience, on figure \ref{fig:Sakiadis}.
This solution was computed by rescaling, with the condition  ${\eta^*_\infty} < {\eta}_{\infty}$, where the chosen truncated boundary was ${\eta^*_\infty} = 10$ in our case.

\subsection{Newton's root-finder}
The same ITM can be applied by using the Newton's root-finder.
This requires a more complex treatment involving a system of six differential equations. 
Let us introduce the auxiliary variables $u_j(\eta)$ for $j = 1, 2, \dots , 6$ defined by
\begin{align}\label{eq:var}
u_1 &= f \ , \quad u_2 = \frac{df}{d\eta} \ , \quad u_3 = \frac{d^2f}{d\eta^2} \ , \nonumber \\[-1ex]
& \\[-1ex]
u_4 &= \frac{\partial u_1}{\partial h^*} \ , \quad u_5 = \frac{\partial u_2}{\partial h^*} \ , \quad u_6 = \frac{\partial u_3}{\partial h^*} \ . \nonumber 
\end{align}
Now, the related IVP is given by
\begin{align}\label{eq:IVP6}
& \frac{du_1^*}{d\eta^*} = u_2^* \ , \nonumber \\
& \frac{du_2^*}{d\eta^*} = u_3^* \ , \nonumber \\
& \frac{du_3^*}{d\eta^*} = -\frac{1}{2}u_1^* u_3^* \ , \nonumber \\[-1.5ex]
&\\[-1.5ex]
& \frac{du_4^*}{d\eta^*} = u_5^* \ , \nonumber \\
& \frac{du_5^*}{d\eta^*} = u_6^* \ , \nonumber \\
& \frac{du_6^*}{d\eta^*} = -\frac{1}{2} \left(u_4^* u_3^* + u_1^* u_6^* \right) \ , \nonumber 
\end{align}
%with initial conditions: 
$u_1^*(0) = 0 \ , \ u_2^*(0) = h^{*1/2} \ , \ u_3^*(0) = -1 \ , \ 
u_4^*(0) = 0 \ , \ u_5^*(0) = \frac{1}{2} h^{*-1/2} \ , \ u_6^*(0) = 0$.

In order to apply the Newton's root-finder, at each iteration, we have to compute the derivative
with respect to $h^*$ of the transformation function $\Gamma$.
In our case, replacing equation (\ref{eq:lambda}) into (\ref{eq:AA2.9}), the transformation function is given by 
\begin{equation}\label{eq:Gamma}
\Gamma(h^*) = \left[u_2^*(\eta_\infty^*)+h^{*1/2}\right]^{-2} h^* -1 \ ,
\end{equation}
and its first derivative can be easily computed as
\begin{equation}\label{eq:dGamma}
\frac{d\Gamma}{dh^*}(h^*) = \left[u_2^*(\eta_\infty^*)+h^{*1/2}\right]^{-2}\left\{1-2 \left[u_5^*(\eta_\infty^*)+\frac{1}{2}h^{*-1/2}\right]\left[u_2^*(\eta_\infty^*)+h^{*1/2}\right]^{-1} h^*\right\}  \ .
\end{equation}
The convergence criterion is again given by (\ref{eq:conv}).
In table \ref{tab:itera2} we list the iterations of our ITM.
\begin{table}[!htb]
\caption{Iterations and numerical results: Newton's root-finder.}
\begin{center}{ \renewcommand\arraystretch{1.3}
%\begin{tabular}{ccccr@{.}lr @{.} lcc} 
\begin{tabular}{rr@{.}lcr@{.}lr@{.}l} 
\toprule%
{$j$}
&
\multicolumn{2}{l}%
{$h_j^*$} 
& {$\lambda_j$}
& \multicolumn{2}{c}%
{$\Gamma(h_j^*)$} 
& \multicolumn{2}{c}%
{${\displaystyle \frac{d^2f}{{d\eta}^2}} (0)$} \\[1.5ex]
\toprule%
% & & & & & \\
$0$  & $2$ &$5$           & $1.061732$ & $0$  & $967345$ & $-0$ & $835517$ \\
$1$  & $2$ & $634888$ & $1.166846$ & $0$ & $421371$ & $-0$ & $629447$ \\
$2$ & $2$ & $812401$ & $1.255130$ & $0$ & $133241$ & $-0$& $505747$\\
$3$ & $2$ & $929233$ & $1.301740$ & $0$ & $020134$ & $-0$ & $453344$ \\
$4$ & $2$ & $953635$ & $1.310767$ & $5$ & $88\mbox{D}{-04}$ & $-0$ & $444042$ \\
$5$ & $2$ & $954391$ & $1.311043$ & $5$ & $26\mbox{D}{-07}$ & $-0$ & $443761$ \\
$6$ & $2$ & $954391$ & $1.311043$ & $4$ & $36\mbox{D}{-13}$ & $-0$ & $443761$ \\
\bottomrule
\end{tabular}}
\label{tab:itera2}
\end{center}
\end{table}

As can be easily seen we get the same numerical results already obtained by the secant method but with a smaller numbers of iterations, that is 5 iterations in table \ref{tab:itera2} compared to 8 iterations in table \ref{tab:itera}.

\section{Conclusions}
The applicability of a non-ITM to the Blasius problem is a consequence of the invariance of the governing differential equation and initial conditions with respect to a scaling group and the non-invariance of the asymptotic boundary condition.
Several problems in boundary-layer theory lack this kind of invariance plus non-invariance and cannot be solved by non-ITMs. 
To overcome this drawback, we can modify the problem at hand by introducing a numerical parameter $ h $, and require the invariance of the modified problem with respect to an extended scaling transformation involving $ h $, see \cite{Fazio:1996:NAN,Fazio:1997:NTE} for the application of this idea to classes of problems.
Here we show how this ITM can be applied to problems with a homogeneous boundary conditions at infinity. 
Moreover, we indicate how to couple our method with the Newton's root-finder.
As far as the choice of a root-finder for the ITM is concerned, we may notice that, if we limit ourselves to consider a scalar nonlinear function, then the secant method, that use one function evaluation per iteration, has an efficiency index higher than the Newton method, where we need at each iteration two function evaluations, as reported by Gautschi \cite[pp. 225-234]{Gautschi:1997:NAI}.
On the other hand, the Newton method can be preferable since it requires only one initial guess.
If we apply these methods to the solution of BVPs, jointly with a shooting or an ITM, then at each iteration the computational cost is by far higher than one or two function evaluations, and as a consequence the Newton's method might be more efficient than the secant one.
\begin{table}[!htb]
\caption{Comparison of the velocity gradient at the plate %at the wall 
and truncated boundary $\eta_\infty $ for the Sakiadis problem.}
\begin{center}{ \renewcommand\arraystretch{1.3}
%\begin{tabular}{ccccr@{.}lr @{.} lcc} 
\begin{tabular}{cccccccc} 
\toprule%
\multicolumn{2}{c}%
{Sakiadis \cite{Sakiadis:1961:BLBa}} 
& \multicolumn{2}{c}%
{Ishak et al. \cite{Ishak:2007:BLM}} 
&\multicolumn{2}{c}%
{Cortell \cite{Cortell:2010:NCB}}
& \multicolumn{2}{c}%
{} \\ 
 \multicolumn{2}{c}%
{} & \multicolumn{2}{c}%
{Finite Difference$^*$} 
&\multicolumn{2}{c}%
{Simple Shooting} 
& \multicolumn{2}{c}%
{ITM} \\ 
\toprule%
{$ \eta_\infty $} &
{$ {\displaystyle \frac{d^2 f}{d\eta^2}} (0) $} &
{$ \eta_\infty $} &
{$ {\displaystyle \frac{d^2 f}{d\eta^2}} (0) $} & %\multicolumn{2}{c}%
{$ \eta_\infty $} & %\multicolumn{2}{c}%
{$ {\displaystyle \frac{d^2 f}{d\eta^2}} (0) $} & 
{$ \eta^*_\infty $} &
{$ {\displaystyle \frac{d^2 f}{d\eta^2}} (0) $} \\[1.5ex] \hline
 &  %& 0 
$-0.44375$ &  & $-0.4438$ & $20$ & $-0.443747$ & $10$ &  $-0.443761$\\
\bottomrule
\end{tabular}}
\label{tab:comp}
\end{center}
\noindent
$^*$ Keller's second order box scheme \cite{Keller74}.
\end{table}
In table \ref{tab:comp} we propose a comparison between our results and those reported in literature.
Our numerical results compare well with those obtained by other authors.

For the Sakiadis flow, which is obtained in the study of a moving flat plate, the modulus of the wall shear $\left|\frac{d^2f}{d\eta^2}(0)\right| = 0.443761$ is larger in comparison to the Blasius flow, for a static flat plate, where we have found $\frac{d^2f}{d\eta^2}(0) = 0.332057$, see \cite{Fazio:1992:BPF}.
The increase of the wall shear can be easily computed by
\[
33.64\% = \frac{|0.332057-0.443761|}{
|0.332057|} 100 \ .
\]
This trend was predicted by Sakiadis \cite{Sakiadis:1961:BLBa} theoretically. 
He proved an increase of about 34\% in the wall shear, see also Sadeghy and Sharifi \cite{Sadeghy:2004:LSS} or Cortell \cite{Cortell:2010:NCB}.
This result was confirmed by Tsou and Goldstein \cite{Tsou:1967:FHT} experimentally.

\bigskip
\bigskip

\noindent
{\bf Acknowledgements.} This work was supported by the University of Messina.

%\clearpage
%\newpage

\end{document}